\documentclass[11pt]{article}
\usepackage[margin=1in]{geometry}
\usepackage{amsmath,amssymb,amsthm,mathtools,amsfonts}
\usepackage[colorlinks=true,linkcolor=blue,citecolor=blue,urlcolor=blue]{hyperref}
\usepackage{authblk}
\usepackage[sort]{cite}
\usepackage{tikz}
\usetikzlibrary{backgrounds,fit,positioning}
\usepackage{graphicx}
\usepackage{amsmath}
\usepackage{mathrsfs}
\usepackage{empheq}

\usepackage{pgfplots}
\pgfplotsset{compat=1.17}

\allowdisplaybreaks
\numberwithin{equation}{section}

\newtheorem{theorem}{Theorem}
\newtheorem{lemma}[theorem]{Lemma}

\newtheorem{problem}[theorem]{Problem}

\newtheorem{conjecture}
{Conjecture}

\theoremstyle{remark}

\newcommand \aln[2]
{
\begin{align}\label{#1}
#2
\end{align}
}

\newcommand \red[1] {{\color{red}{#1}}}

\newcommand{\I}{\mathcal I}

\newcommand{\N}{\mathbb N}
 \newcommand{\R}{\mathbb R}

\usepackage{nccmath}  

\begin{document}
	
	\title{A sharp form of Welsh's log-concavity inequality for chromatic polynomials}
	
\author[1\thanks{Corresponding author. Email:
			meiqiaozhang95@163.com and meiqiaozhang@xmu.edu.cn.}]{Meiqiao Zhang}
		
\author[2\thanks{Email: fengming.dong@nie.edu.sg 
		and donggraph@163.com.}
	]{Fengming Dong}
	
	\affil[1]{\small School of Mathematical Sciences, Xiamen  University, China}

	\affil[2]{\small
		National Institute of Education,
		Nanyang Technological University, 
		Singapore}

	\date{}
	
	\maketitle

\begin{abstract}
Let $G$ be a simple graph of order $n\ge 1$, and let $P_G(x)$ be its chromatic polynomial. In the early 1970s, Welsh conjectured that $P_G(k)^2 \ge P_G(k-1)P_G(k+1)$ holds for any positive integer $k$. This conjecture was also independently proposed by Brenti in 1992, while in 1997, Seymour disproved it for $k=6$ by constructing an infinite family of graphs of sufficiently large order as counterexamples. In 2005, Dong, Koh and Teo conjectured that $P_G(x)^2 \ge P_G(x-1)P_G(x+1)$ holds for all real $x\ge n$. In this paper, we prove that the inequality holds for all real  $x\ge n-1$. Equality holds if and only if $n\ge 2$, $G$ is complete, and $x=n-1$. Moreover,  the bound $x\ge n-1$ is best possible for every $n\ge 2$.
\end{abstract}

\smallskip
\noindent \textbf{Keywords:} chromatic polynomial, log-concavity,  graph coloring

\smallskip
\noindent \textbf{Mathematics Subject Classification: 05C31, 05C15}

%\tableofcontents 

\section{Introduction}\label{sec:intro}
In this article, we consider simple graphs only.
For any graph $G$, denote by $V(G)$ and $E(G)$ the vertex set and the edge set of $G$, respectively.
For any $v\in V(G)$, let \(N_G(v)\) denote the set of neighbors of \(v\) in \(G\).
For any $S\subseteq V(G)$, write $G[S]$ for the subgraph of $G$ induced by $S$, $G-S$ for $G[V(G)\setminus S]$, and 
$G_S$ for the graph obtained from $G$ by adding a new vertex 
$v$ and new edges joining 
$v$ to all vertices in $S$.
Let $N$ denote the graph without vertices and $K_n$ denote the complete graph of order $n$. 
Denote by $\N$ and $\R$ the sets of positive integers and real numbers, respectively. For any $k\in\N$ and $x\in \R$, let $[k]=\{1,\dots,k\}$, and  let $(x)_k=x(x-1)\cdots(x-k+1)$ be the $k$-th falling factorial polynomial. Set 
$[0]=\emptyset$ and
$(x)_0=1$ by convention.

%Denote by $\mathbb N, \mathbb R, \mathbb C$ the sets of positive integers,  real numbers and complex numbers, respectively. For any $m\in\N$, let $[m]=\{1,\ldots,m\}$.

Let $G$ be a graph of order $n\ge 1$. 
For any $k\in \N$, a 
{\it proper $k$-coloring} of $G$ is a mapping $\theta:V(G)\rightarrow [k]$ such that $\theta(u)\neq \theta(v)$ for all $uv\in E(G)$. 
The \textit{chromatic polynomial} $P_G(x)$ of $G$ is defined to be the unique polynomial that counts the number of proper $k$-colorings of $G$ for each $x=k \in \N$. Also, let $P_N(x)=1$ by convention. 
The chromatic polynomial was introduced by Birkhoff in~\cite{birk} as a tool for attacking the Four-Color Conjecture, but later became an important object of study due to its elegant  combinatorial and algebraic properties~\cite{birk, birk2, dong1, DKT, Jackson2015, rea1, rea2}.

In particular, the study on log-concavity and unimodality of chromatic polynomials has a long history.
For any $\ell\in\N$, a finite sequence $(a_0,a_1,\ldots,a_\ell)$ of nonnegative
numbers is said to be \textit{log-concave} if
$
a_k^2\ge a_{k-1}a_{k+1}$ for all $1\le k\le \ell-1$, and is \textit{unimodal} if there exists $0\le k\le\ell$ such that
$a_0\le \cdots\le a_k\ge\cdots\ge a_\ell.$ 
Obviously, a log-concave sequence with no internal zeros is unimodal. 
In 1968, Read~\cite{rea1} conjectured that the absolute values of the coefficients of the chromatic polynomial of every graph form a unimodal sequence, while the stronger log-concavity property was subsequently established by Huh~\cite{Huh} in 2012 using methods from algebraic geometry.
Similar research also appears naturally for the coefficient sequences obtained by expanding the chromatic polynomial in other bases~\cite{Brenti,Brown,Zhang}.

Meanwhile, the log-concavity of the values of chromatic polynomials has been extensively studied. In the early 1970s, Welsh proposed the following conjecture concerning the log-concavity of evaluations at positive integers.

\begin{conjecture}\label{conj}
For any graph $G$ and $k\in \N$, $P_G(k)^2 \ge P_G(k-1)P_G(k+1)$. 
\end{conjecture}
Brenti~\cite{Brenti} independently proposed the same conjecture in 1992 and verified it for a broad class of graphs.
However, in 1997, Seymour~\cite{Seymour} disproved Conjecture~\ref{conj} by constructing graphs $G$ of arbitrarily large order for which $P_G(6)^2<P_G(5)P_G(7)$.
Since these counterexamples violate the conjecture at values of $k$ that are small relative to the orders of graphs, it is natural to ask whether the conjecture holds for all sufficiently large integers $k$.
Further, Dong, Koh and Teo~\cite{DKT} proposed the following conjecture.

\begin{conjecture}[\cite{DKT}]
\label{conj1}
For any graph $G$ of order $n$, $P_G(x)^2 \ge P_G(x-1)P_G(x+1)$ for all real $x\ge n$. 
\end{conjecture}

To support Conjecture~\ref{conj1}, Dong, Koh and Teo~\cite{DKT} proved that the inequality holds for all real $x\ge\max\{|V(G)|-1,\sqrt{2}(|E(G)|-|V(G)|+\frac52)\}$ when $G$ is connected. Subsequently, 
Fadnavis~\cite{Fadnavis} established the inequality for every graph $G$ and all integers $x\ge C\Delta(G)+1$, where $\Delta(G)$ is the maximum degree of $G$ and $0<C<10$ is a constant.

In this paper, we confirm Conjecture~\ref{conj1} by the following stronger result.

\begin{theorem}\label{main}
Let $G$ be a graph of order $n\ge 1$.
For every real $x\ge n-1$,
\aln{impor}{
  P_G(x)^2\ge P_G(x-1)P_G(x+1),
}
where the equality holds if and only if $G\cong K_n$ with $n\ge 2$ and $x=n-1$.
\end{theorem}

\noindent{\bf Remark}.
Note that for $n\ge 2$,
\aln{complete}
{
P_{K_n}(x)^2- P_{K_n}(x-1)P_{K_n}(x+1)
%& =(x)_n^2-(x-1)_n(x+1)_n \nonumber\\
=nx(x-n+1)\left ((x-1)_{n-2}\right )^2.
%n(x)_{n-1}(x-1)_{n-1},
}
By \eqref{complete}, 
(\ref{impor}) clearly holds for $G=K_n$ whenever $x\ge n-1$, with equality at $n\ge 2$ and $x=n-1$. 
Moreover, 
\eqref{complete} shows that 
(\ref{impor}) fails for every $n\ge 2$ and $n-2<x<n-1$.
 Thus the condition $x\ge n-1$ in Theorem~\ref{main} is the best possible.

The proof of Theorem~\ref{main} will be given in Section~\ref{sec3} after introducing some preliminary results in Section~\ref{sec2} and
two key lemmas in Section~\ref{key-lem}.
Further problems
are proposed in Section~\ref{further}.

\section{Preliminaries}\label{sec2}

In this section, we establish a series of identities involving chromatic polynomials.

For any graph $G$ of 
order  $n\ge 1$ 
and each $k\in [n]$,
%integer $1\le k\le n$, 
let $\alpha_k(G)$ denote the number of partitions of $V(G)$ into $k$ nonempty independent sets. Since each proper $k$-coloring of $G$ induces a partition of $V(G)$ into at most $k$ independent sets, 
the chromatic polynomial has the following well-known
falling-factorial expansion
(see, e.g., \cite{Brown, DKT, rea1}):
%the chromatic polynomial has the following well-known expansion:
\aln{expre1}
{
 P_G(x)=\sum_{k=1}^n \alpha_k(G)(x)_k.
}
It is clear that $\alpha_n(G)=1$,
and $\alpha_{n-1}(G)\ge 1$ whenever $G\not\cong K_n$.
%$P_G(x)\ge (x)_n$,  and $P_G(x)\ge (x)_n+(x)_{n-1}$ if $G\not \cong K_n$. 
Therefore,  
\aln{positive}
{P_G(x)>0~~\text{whenever}~x> n-1,~\text{or}~x>n-2~\text{and}~G\not \cong K_n.
}

Further, for any graph $G$ of order $n\ge 0$,
let $\I(G)$ be the set of independent sets of $G$. By convention, $\emptyset\in \I(G)$.
For any $x\in \R$, $I\in \I(G)$
and $u\in V(G)$, define  
\begin{empheq}[left=\empheqlbrace]{align}
	D_G(x)
	&= P_G(x)^2-P_G(x-1)P_G(x+1).
	\label{eq:D}\\
	\phi_G(x,I)
	&= P_G(x)P_{G-I}(x-1)
	-P_G(x-1)P_{G-I}(x).
	\label{eq:phi}\\
	\ell_{G,u}(x)
	&= \sum_{I\in\I(G)\atop u\in I}\phi_G(x,I).
	\label{eq:ell}\\
	\ell_G(x)&=\frac{P_G(x)^2}{x}
	-\frac{P_G(x-1)P_G(x+1)}{x+1}\quad\text{for}~ x>0.
		\label{eq:ell-0}
\end{empheq}
Then we have the following lemma. 

\begin{lemma}\label{key1}
	Let $G$ be a graph of order $n\ge 1$, let $u\in V(G)$,
	and let $S\subseteq V(G)$. Then for any $x\in \R$, the following identities hold:
	\begin{fleqn}[0pt]
	%	\begin{alignat}{2}
	\begin{alignat}{2}
		&\textnormal{(i)}\quad
		& P_G(x)
		&= \sum_{I\in\mathcal I(G)} P_{G-I}(x-1).
		\notag\\[3pt]
		&\textnormal{(ii)}\quad
		& P_G(x)
		&= x\sum_{\substack{I\in\mathcal I(G)\\ u\in I}}
		P_{G-I}(x-1).
		\notag\\[3pt]
		&\textnormal{(iii)}\quad
		& P_{G_S}(x)
		&= x\sum_{\substack{I\in\mathcal I(G)\\ I\cap S=\varnothing}}
		P_{G-I}(x-1).
		\notag\\[3pt]
		&\textnormal{(iv)}\quad
		& D_G(x)
		&= \sum_{I\in\mathcal I(G)} \phi_G(x,I).
		\notag\\[3pt]
		&\textnormal{(v)}\quad
		& \ell_{G,u}(x)
		&=\frac{P_G(x)^2}{x}
		-\frac{P_G(x-1)P_G(x+1)}{x+1}=\ell_G(x)
	~\text{when}~ x>0.
		\notag
	\end{alignat}
\end{fleqn}
\end{lemma}

\begin{proof}
Since every proper $k$-coloring of $G$ is uniquely determined by the independent set of vertices receiving the first color and a proper $(k-1)$-coloring of the remaining graph, (i) is obvious. 

Similarly, for each $u\in V(G)$, every proper $k$-coloring of $G$ is uniquely determined by the choice of a color for $u$, the independent set of vertices receiving that color, and a proper $(k-1)$-coloring of the remaining graph. (ii) follows immediately.

For (iii), as $G_S$ is the graph obtained from $G$ by adding a new vertex $v$ and attaching $v$ to each vertex in $S$,
applying (ii) to $G_S$ with $u=v$ gives 
$$
P_{G_S}(x)= x\sum_{I'\in \I(G_S)\atop v\in I'}
P_{G_S-I'}(x-1)=x\sum_{I\in\I (G)\atop I\cap S=\emptyset}P_{G-I}(x-1).
$$

To establish (iv), 
we use (\ref{eq:D}), (i) and (\ref{eq:phi}) to obtain that
\aln{}{
 D_G(x) & =  P_G(x)^2-P_G(x-1)P_G(x+1)\nonumber\\
& =   P_G(x)\sum_{I\in \I(G)}  P_{G-I}(x-1)-P_G(x-1)\sum_{I\in \I(G)}  P_{G-I}(x)\nonumber\\
& =   \sum_{I\in\I(G)}  \phi_G(x,I).\nonumber
}

Finally, for (v), we combine  (\ref{eq:ell}), (\ref{eq:phi}), (ii) and (\ref{eq:ell-0}) to obtain that
\aln{}{
\ell_{G,u}(x)
  & = P_G(x)\sum_{I\in\I(G)\atop u\in I}P_{G-I}(x-1)-P_G(x-1)\sum_{I\in\I(G)\atop u\in I}P_{G-I}(x)\nonumber\\
  & = \frac{P_G(x)^2}{x}
     -\frac{P_G(x-1)P_G(x+1)}{x+1}\nonumber\\
&=\ell_G(x),
     \nonumber
}
where the second equality holds as $x>0$. 
\end{proof}

%\red{By (\ref{eq:ell-for}), the value of  $\ell_{G,u}(x)$ is independent of the choice of $u$.}

 By (\ref{eq:D}), proving (\ref{impor}) is equivalent to showing that $D_G(x)\ge 0$. Moreover, since Lemma~\ref{key1}(iv) expresses $D_G(x)$ in terms of $\phi_G(x,I)$ over all $I\in\I(G)$, we now introduce some auxiliary functions to further study $\phi_G(x,I)$.

%\section{Proof of Theorem~\ref{main}}\label{sec3}
%In this section, we prove Theorem~\ref{main} after establishing a key lemma.

For any graph $G$ of order 
$n\ge 0$
and $x\in\R$ such that  $P_G(x)> 0$, 
define
\aln{eq:R}{
R_G(x)
= \frac{P_G(x-1)}{P_G(x)},
}
and for any $S\subseteq V(G)$, %define 
\aln{eq:A}{
	A_G(x,S)
	= \frac{P_{G_S}(x)}{P_G(x)}.
}
Then we have the next lemma.

\begin{lemma}\label{key2}
	Let $G$ be a graph of order $n\ge 1$ and $x\in \R$. 
	The following results hold:
		\begin{enumerate}
	\item  For any $I\in \I(G)$, if  $P_G(x), P_{G-I}(x)>0$,  then 
	\aln{eq10}
	{
	\phi_G(x,I)
=P_G(x)P_{G-I}(x)\left(R_{G-I}(x)-R_G(x)\right).
}
\item  If $P_G(x), P_G(x-1)>0$, 
and $P_{G_S}(x)>0$ for $S\subseteq V(G)$,  then 
	\aln{eq:proA}
	{
R_{G_S}(x)A_G(x,S)
	= R_G(x)A_G(x-1,S).
}
\end{enumerate}
\end{lemma}

\begin{proof}
For (i), suppose $I\in\I(G)$ with $P_G(x), P_{G-I}(x)>0$. Then (\ref{eq:phi}) yields
\begin{align*}
\phi_G(x,I) & =P_G(x)P_{G-I}(x-1)-P_G(x-1)P_{G-I}(x)\nonumber\\
& = P_G(x)P_{G-I}(x) \frac{P_{G-I}(x-1)}{P_{G-I}(x)}-P_{G-I}(x)P_G(x)\frac{P_{G}(x-1)}{P_{G}(x)}\nonumber\\
&  = P_G(x)P_{G-I}(x)\left(R_{G-I}(x)-R_G(x)\right).
\end{align*}
\iffalse 
Moreover, for any $S\subseteq V(G)$, define 
\aln{eq:A}{
A_G(x,S)
= \frac{P_{G_S}(x)}{P_G(x)}.}
\fi

For (ii), suppose $P_G(x), P_G(x-1), P_{G_S}(x)>0$.
Then by (\ref{eq:R}) and (\ref{eq:A}),
$$
R_{G_S}(x)A_G(x,S)= \frac{P_{G_S}(x-1)}{P_G(x)}=R_G(x)A_G(x-1,S).$$
\end{proof}
The identities established
in Lemmas~\ref{key1} and~\ref{key2} 
will be used  in the next two sections.

\section{Key lemmas}
\label{key-lem}

In this section, we  present two key lemmas for proving Theorem~\ref{main} in the next section.

\begin{lemma}
	\label{lem:delete1}
	Let $G$ be a graph of order 
	$n\ge 0$ and 
	$S\subseteq V(G)$.
 Then
 $A_G(x,S)\ge x-|S|$ for
 all $x\ge n$.
\end{lemma}

\begin{proof}
%	We shall prove by induction on $n$. 	
	If $n=0$, then $S=\emptyset$,  $P_G(x)=1$ and $G_S\cong K_1$. Thus $P_{G_S}(x)=x$, 
		implying that 
		$A_G(x,S)=x=x-|S|$.
		Hence the result holds  in this case.
	
	Now suppose %$n\ge1$.
%	If $G\cong K_n$, then $P_G(x)=(x)_n$ and $P_{G_S}(x)=(x-|S|)(x)_n$. Thus, $A_G(x,S)=(x-|S|)$and the result holds  in this case.
%	In the following, assume that %$G\not\cong K_n$ and 
$x\ge n\ge 1$ and $S$ is a fixed 
subset of $V(G)$.
%Then $G_S\not\cong K_{n+1}$.
	%It suffices to show that $A_G(x,S)\ge x-|S|$ and $R_{G_S}(x)\le R_G(x)$. 
	Note that $P_G(x)>0$ %P_{G_S}(x)>0$ 
	by (\ref{positive}), 
	implying that $A_G(x,S)$ %R_{G_S}(x), R_G(x)$ 
	is  well-defined.

	Let %$f:\I(G)\to \R$  with 
	$f(I)=x P_{G-I}(x-1)/P_G(x)$ for each $I\in \I(G)$. 
Since $P_{G-I}(x-1)\ge 0$ by (\ref{expre1}), $f(I)$ is always nonnegative. Moreover, 
due to Lemma~\ref{key1}(i)--(iii) and (\ref{eq:A}),
	\begin{align}
		\sum_{I\in \I(G)}f(I) & = \sum_{I\in \I(G)} \frac{xP_{G-I}(x-1)}{P_G(x)} = \frac{xP_G(x)}{P_G(x)}=x.\label{2-8}\\
\forall u\in V(G),\quad		\sum_{I\in \I(G)\atop u\in I}f(I) & = \sum_{I\in \I(G)\atop u\in I} \frac{xP_{G-I}(x-1)}{P_G(x)} = \frac{P_G(x)}{P_G(x)}=1.
\label{2-10}\\
		\sum_{I\in \I(G)\atop I\cap S=\emptyset}f(I)  & = \sum_{I\in \I(G)\atop I\cap S=\emptyset} \frac{xP_{G-I}(x-1)}{P_G(x)} = \frac{P_{G_S}(x)}{P_G(x)}=A_G(x,S).\label{2-9}
	\end{align}
	%where $u\in V(G)$.
	
	Further,
%	for given $S\subseteq V(G)$,
	 (\ref{2-10}) gives
	\aln{2-11}{
		|S|  = \sum_{u\in S}1 
		= \sum_{u\in S}\sum_{I\in\mathcal I(G)\atop u\in I} f(I) 
		= \sum_{I\in\mathcal I(G)} |I\cap S|\,f(I) 
		\ge \sum_{I\in\mathcal I(G)\atop I\cap S\ne\emptyset} f(I),
	}
	where the last inequality holds as $f(I)\ge 0$ for  each $I\in\I(G)$.
	
	Then by (\ref{2-9}), (\ref{2-8}) and (\ref{2-11}),
	$$
	A_G(x,S) = \sum_{I\in \I(G)\atop I\cap S=\emptyset}f(I)
	=  \sum_{I\in \I(G)}f(I)-\sum_{I\in \I(G)\atop I\cap S\neq \emptyset}f(I)
	\ge x-|S|.
	$$
Hence the result holds.
\end{proof}

\begin{lemma}
	\label{lem:delete}
Let $G$ be a graph of order 
$n\ge 0$ and 
$S\subseteq V(G)$.
Then $R_{G_S}(x)\le R_G(x)$
for every real $x> n$.
If \(G_S\not\cong K_{n+1}\), the inequality also holds at \(x=n\).
%if $G_S\cong K_{n+1}$,and for all real $x\ge n$ otherwise. 
\iffalse 
\begin{enumerate}
\item[(i)]  $A_G(x,S)\ge x-|S|$, and 
%for every real $x\ge n$. 
\item[(ii)] 
%For every real $x>n$, $R_{G_S}(x)\le R_G(x).$
 $R_{G_S}(x)\le R_G(x)$
 when $G_S\not\cong K_{n+1}$, 
 and for $x>n$ otherwise.

%\item[(iii)] If $G_S\not\cong K_{n+1}$,  then $R_{G_S}(n)\le R_G(n)$.
\end{enumerate}\fi
\end{lemma}

\begin{proof}
 We shall prove by induction on $n$. 
 
 If $n=0$, then $S=\emptyset$,  $P_G(x)=1$ and 
 $G_S\cong K_1$. 
Thus $P_{G_S}(x)=x$. 
By \eqref{eq:R}, 
$$
R_{G_S}(x)=(x-1)/x<1=R_G(x),
$$
when $x>n=0$.
Hence the lemma is true in this case.
 
Now suppose $n\ge1$.
If $G\cong K_n$, then $P_G(x)=(x)_n$ and $P_{G_S}(x)=(x-|S|)(x)_n$. 
Thus,  by \eqref{eq:R}, 
$$
R_{G_S}(x)
=\frac{(x-n)(x-1-|S|)}{x(x-|S|)}
\le \frac{x-n}{x}=R_G(x),
$$
whenever $x>n$, or $x=n$ and $G_S\not\cong K_{n+1}$.
Hence,  the lemma holds in this case.

In the following, assume that $G\not\cong K_n$, $S$ is a subset of $V(G)$ and $x\ge n\ge 2$. Then $G_S\not\cong K_{n+1}$. It suffices to prove that  $R_{G_S}(x)\le R_G(x)$.
Note that $P_G(x),  P_{G_S}(x)>0$ follows from (\ref{positive}), which implies that 
both $ R_{G_S}(x)$ and 
$R_G(x)$ are well-defined. 
We complete the proof by two claims below. 
 
 \iffalse 
 \red{
\noindent{\bf Claim 1}.
 $A_G(x,S)\ge x-|S|$. 
}
 
 Let $f:\I(G)\to \R$ 
 with $f(I)=x P_{G-I}(x-1)/P_G(x)$ for each $I\in \I(G)$. 
 Since $P_{G-I}(x-1)>0$ follows from (\ref{positive}), $f(I)$ is always positive. Moreover, by Lemma~\ref{key1}(i)--(iii) and (\ref{eq:A}),
\begin{align}
\sum_{I\in \I(G)}f(I) & = \sum_{I\in \I(G)} \frac{xP_{G-I}(x-1)}{P_G(x)} = \frac{xP_G(x)}{P_G(x)}=x,\label{2-8}\\
\sum_{I\in \I(G)\atop I\cap S=\emptyset}f(I)  & = \sum_{I\in \I(G)\atop I\cap S=\emptyset} \frac{xP_{G-I}(x-1)}{P_G(x)} = \frac{P_{G_S}(x)}{P_G(x)}=A_G(x,S),\label{2-9}\\
\sum_{I\in \I(G)\atop u\in I}f(I) & = \sum_{I\in \I(G)\atop u\in I} \frac{xP_{G-I}(x-1)}{P_G(x)} = \frac{P_G(x)}{P_G(x)}=1,\label{2-10}
\end{align}
where $u\in V(G)$.
Further, (\ref{2-10}) gives
\aln{2-11}{
|S|  = \sum_{u\in S}1 
  = \sum_{u\in S}\sum_{I\in\mathcal I(G)\atop u\in I} f(I) 
  = \sum_{I\in\mathcal I(G)} |I\cap S|\,f(I) 
  \ge \sum_{I\in\mathcal I(G)\atop I\cap S\ne\emptyset} f(I).
}
Then by (\ref{2-9}), (\ref{2-8}) and (\ref{2-11}),
$$
A_G(x,S) = \sum_{I\in \I(G)\atop I\cap S=\emptyset}f(I)
=  \sum_{I\in \I(G)}f(I)-\sum_{I\in \I(G)\atop I\cap S\neq \emptyset}f(I)
  \ge x-|S|.
$$
\red{Hence Claim 1 holds.}
\fi

\noindent{\bf Claim 1}.
$\phi_G(x-1, I)\ge 0$ for 
each $I\in \I(G)$.

%Now it remains to prove $R_{G_S}(x)\le R_G(x)$. 

%We shall show that $\phi_G(x-1, I)\ge 0$ for all $I\in \I(G)$ first.

By (\ref{eq:phi}),
$\phi_G(x-1,\emptyset)=0$.
Now consider any nonempty set  $I=\{v_1,\dots,v_q\}$ 
	in $\I(G)$.
Let $G_0=G$ and $G_i=G-\{v_1,\dots,v_i\}$ for each $i\in [q]$. 
For each $i\in [q]$, $G_i$ has order $n-i$ and $(G_i)_{N_G(v_i)}\cong G_{i-1}$. 
Then, $x-1\ge n-1=|V(G_1)|$
and $G_0=G\not\cong K_n$,
and $x-1\ge n-1>|V(G_i)|$
if $i\ge 2$.
For each $i\in [q]$, 
by the induction hypothesis on $G_i$ with $S=N_G(v_i)$, 
we have 
$R_{G_{i-1}}(x-1)\le R_{G_i}(x-1)$.
 Hence  
\aln{chain}
{R_G(x-1)=R_{G_0}(x-1)
	\le R_{G_1}(x-1)\le \cdots \le  R_{G-I}(x-1).}
Recall that $P_G(x-1), P_{G-I}(x-1)>0$ by \eqref{positive}. 
Thus (\ref{eq10}) and (\ref{chain}) together yield
$$
 \phi_G(x-1,I) = P_G(x-1) P_{G-I}(x-1)\bigl(R_{G-I}(x-1)-R_G(x-1)\bigr)\ge  0.
$$

Thus Claim 1 holds.

Claim 1 implies that $\ell_{G,u}(x-1)\ge 0$ for 
each $u\in V(G)$ by (\ref{eq:ell}). Then $\ell_{G}(x-1)\ge 0$ follows from  Lemma~\ref{key1}(v).

	\noindent{\bf Claim 2}.
	$R_{G_S}(x)\le R_G(x)$.

By (\ref{eq:phi}), Lemma~\ref{key1}(iii), (\ref{eq:A}) and (\ref{eq:ell-0}), we have
\aln{ieq16}{
	\sum_{I\in \I(G)\atop I\cap S=\emptyset} \phi_G(x-1,I)
	&= P_G(x-1)\sum_{I\in \I(G)\atop I\cap S=\emptyset} P_{G-I}(x-2)
	- P_G(x-2)\sum_{I\in \I(G)\atop I\cap S=\emptyset} P_{G-I}(x-1)\nonumber\\
	&= \frac{P_G(x-1)P_{G_S}(x-1)}{x-1}
	- \frac{P_G(x-2)P_{G_S}(x)}{x}\nonumber\\
	&= \frac{P_G(x-1)^2}{x-1}A_G(x-1,S)
	- \frac{P_G(x-2)P_G(x)}{x}A_G(x,S)\nonumber\\
	&= \frac{P_G(x-1)^2}{x-1}A_G(x-1,S)
	+ \left(\ell_{G}(x-1)-\frac{P_G(x-1)^2}{x-1}\right)A_G(x,S)\nonumber\\ 
	&= \ell_{G}(x-1) A_G(x,S)
	+ \frac{P_G(x-1)^2}{x-1}\bigl(A_G(x-1,S)-A_G(x,S)\bigr),
}
where the second and fourth equalities are valid as $x\ge n\ge 2$, and the third equality holds as $P_G(x-1), P_G(x)>0$.

Also, the fact that $\phi_G(x-1,\emptyset)=0$, Claim 1, (\ref{eq:ell}) and Lemma~\ref{key1}(v) yield that
\aln{ieq15}{
  \sum_{I\in \I(G)\atop I\cap S=\emptyset}  \phi_G(x-1,I)
 & \le \sum_{I\in \I(G)\atop I \cap S=\emptyset} |I|\,\phi_G(x-1,I)\nonumber\\
 & = \sum_{u\in V(G)\setminus S}
       \sum_{I\in \I(G)\atop u\in I, I\cap S=\emptyset} \phi_G(x-1,I)\nonumber\\
 & \le \sum_{u\in V(G)\setminus S}\sum_{I\in \I(G)\atop u\in I} \phi_G(x-1,I)\nonumber\\
&  = \sum_{u\in V(G)\setminus S} \ell_{G,u}(x-1)\nonumber\\
&  = (n-|S|)\ell_{G}(x-1).
}

Then by (\ref{ieq15}) and (\ref{ieq16}),
$$\ell_G(x-1) A_G(x,S)
+ \frac{P_G(x-1)^2}{x-1}\bigl(A_G(x-1,S)-A_G(x,S)\bigr)
\le (n-|S|)\ell_G(x-1),
$$
implying that 
\aln{ine18}{
	\frac{P_G(x-1)^2}{x-1}\bigl(A_G(x,S)-A_G(x-1,S)\bigr)
	&  \ge   \bigl(A_G(x,S)-n+|S|\bigr) \ell_G(x-1)\nonumber\\
	&  \ge  (x-|S|-n+|S|) \ell_G(x-1) \nonumber\\
	&  \ge  0,
}
where the penultimate inequality follows from Lemma~\ref{lem:delete1} and the fact that 
$\ell_G(x-1)\ge 0$.

Since $x-1>0$ and 
$ P_G(x-1)>0$, (\ref{ine18}) indicates $A_G(x,S)\ge A_G(x-1,S)$. Moreover, $P_G(x-1), P_G(x), P_{G_S}(x)>0$ imply $R_G(x), A_G(x,S)>0$, which together with (\ref{eq:proA}) yield
$$
R_{G_S}(x)=\frac{A_G(x-1,S)}{A_G(x,S)}R_G(x)\le R_G(x).
$$
Claim 2 holds.

Hence 
the lemma follows by induction.
\iffalse 
where the last inequality follows from  (\ref{ieq15}).

By Lemma~\ref{key1}(v), %and~\eqref{positive}, 
as $x\ge n$,
\aln{}
{
\frac{P_G(x-1)^2}{x-1}\cdot 
A_G(x,S)
&=
\left (\ell_G(x-1)
+\frac{P_G(x-2)P_G(x)}{x}
\right )\cdot \frac{P_{G_S}(x)}{P_G(x)}
\nonumber
\\
&=\ell_G(x-1)A_G(x,S)
+\frac{P_G(x-2)P_{G_S}(x)}{x}
.
}
and by \eqref{eq:A}, 
$$
\frac{P_G(x-1)^2}{x-1}
A_G(x-1,S)
=\frac{P_G(x-1)P_{G_S}(x-1)}{x-1}
$$

Then by Lemma~\ref{key1}(v),  (i) and (\ref{ieq16}),
\aln{ine18}{
  \frac{P_G(x-1)^2}{x-1}\bigl(A_G(x,S)-A_G(x-1,S)\bigr)
  &=
  \\
&  \ge   \bigl(A_G(x,S)-n+|S|\bigr) \ell_{G,u}(x-1)\nonumber\\
&  \ge  (x-|S|-n+|S|) \ell_{G,u}(x-1) \nonumber\\
&  \ge  0.
}

Since $x-1, P_G(x-1)>0$, (\ref{ine18}) indicates $A_G(x,S)\ge A_G(x-1,S)$. Moreover, $P_G(x-1), P_G(x), P_{G_S}(x)>0$ imply $R_G(x), A_G(x,S)>0$, which together with (\ref{eq:proA}) yield
$$
  R_{G_S}(x)=\frac{A_G(x-1,S)}{A_G(x,S)}R_G(x)\le R_G(x).
$$
Hence (ii) and (iii) hold. The lemma follows by induction.\fi
\end{proof}

\section{Proof of the main result %Theorem~\ref{main}
}
\label{sec3}

Now we prove Theorem~\ref{main}.

\begin{proof}[Proof of Theorem~\ref{main}]
Let $G$ be a graph
of order $n\ge 1$ and $x\ge n-1$. 
Then
it is clear by (\ref{complete}) that (\ref{impor}) holds when $G\cong K_n$, where the equality holds if and only if $n\ge 2$ and $x=n-1$. 

Now assume $G\not\cong K_n$. Then $n\ge 2$.
Denote by $\omega(G)$ the clique number of $G$. Then $1\le \omega(G)\le n-1$. 
In the following, 
we shall prove that
\aln{final}{D_G(x)
  \ge \frac{\omega(G)P_G(x)^2}{x(x+1-\omega(G))}>0.}
Since $x\ge n-1\ge \omega(G)>0$ and $P_G(x)>0$ by (\ref{positive}), it remains to show the first inequality. 

Note that $\phi_G(x,\emptyset)=0$ by (\ref{eq:phi}). 
Also, for each nonempty $I\in\I(G)$, $G$ is the graph obtained from $G-I$ by adding and attaching vertices in $I$ one by one properly, where $|V(G-I)|\le n-1\le x$. Since $G$ is not complete, $R_G(x)\le R_{G-I}(x)$ follows from  Lemma~\ref{lem:delete}. Notice that $P_{G-I}(x)>0$ by (\ref{positive}).
Then (\ref{eq10}) gives
$$
\phi_G(x,I) = P_G(x)P_{G-I}(x)\left(R_{G-I}(x)-R_G(x)\right)\ge 0.
$$
Thus 
	$\ell_{G,u}(x) \ge 0$ 
	for each $u\in V(G)$
	by (\ref{eq:ell}), which further implies that 
	$\ell_G(x)
\ge 0$ 
due to Lemma~\ref{key1}(v).

Let $C$ be a clique in $G$ with $|C|=\omega(G)$. 
Since $|C\cap I|\le 1$ for each 
$I\in \I(G)$, 
Lemma~\ref{key1}(iv)(v) and \eqref{eq:ell} yield
\aln{ieq20}{
  \omega(G)\ell_{G}(x)
  =\sum_{u\in C}\sum_{I\in\I(G)\atop u\in I}\phi_G(x,I)
  =\sum_{I\in\I(G)\atop I\cap C\ne\emptyset}\phi_G(x,I)
  \le  \sum_{I\in\I(G)}\phi_G(x,I)=D_G(x).
}
On the other hand, by Lemma~\ref{key1}(v) and  (\ref{eq:D}),
\aln{ieq6}{
  (x+1)\ell_{G}(x)
  = (x+1)\left(\frac{P_G(x)^2}{x}
    -\frac{P_G(x-1)P_G(x+1)}{x+1}\right) = D_G(x)+\frac{P_G(x)^2}{x}.
}
Consequently, (\ref{ieq20}) and (\ref{ieq6}) together indicate $(x+1-\omega(G))\ell_{G}(x)\ge P_G(x)^2/x$. As
$x+1-\omega(G)>0$, (\ref{ieq20}) further implies that 
$$
  D_G(x)\ge\omega(G)\ell_{G}(x)
  \ge\frac{\omega(G)P_G(x)^2}{(x+1-\omega(G))x}.
$$
This establishes (\ref{final}), and hence (\ref{impor}), completing the proof of Theorem~\ref{main}.
\end{proof}

\iffalse 
\red{
For any integers $n\ge 2$ and 
$0\le d\le n$, 
let $G_{d,n}$ be the graph obtained 
from $K_n$ by adding a new vertex $w$ and $d$ new edges
adding $w$ to $d$ vertices in $K_n$.
Then,
$$
P_{G_{d,n}}(x)=(x-d)(x)_n.
$$
Thus,
\aln{se4-e1}
{
&P_{G_{n-2,n}}(x)^2 -
P_{G_{n-2,n}}(x+1)P_{G_{n-2,n}}(x-1)
\nonumber \\
&=(x-(n-2))^2 ((x)_n)^2
-(x+1-(n-2))(x+1)_n
(x-1-(n-2))(x-1)_n
\nonumber \\
&=((x)_n)^2
\left [ (x-n+2)^2 
-(x-n+3)(x-n+1)\cdot \frac{(x+1)(x-n)}{(x-n+1)x}
\right ]
\nonumber \\
&=((x)_n)^2
\frac{(n+1)x -n^{2}+3 n}{x}.
}
Thus, if $n\ge 4$ and $0<x<\frac{n \left(n-3\right)}{n+1}$, then 
$$
P_{G_{n-2,n}}(x)^2 
<
P_{G_{n-2,n}}(x+1)P_{G_{n-2,n}}(x-1).
$$
Similarly, for $n\ge 3$, 
we can show that 
\aln{se4-e2}
{
	&P_{G_{n-3,n}}(x)^2 -
	P_{G_{n-3,n}}(x+1)
	P_{G_{n-3,n}}(x-1)
	\nonumber \\
&=((x)_n)^2\cdot 
\frac{(n+1)x^{2}-
(2 n^{2}-5 n x-1)x+n^{3}-6 n^{2}+8 n}{x(x-n+1)}.
}
The two roots of the equation 
$
(n+1)x^{2}-
(2 n^{2}-5 n x-1)x+n^{3}-6 n^{2}+8 n$ are 
$$
\frac{2 n^{2}-5 n-1+\sqrt{13 n^{2}-22 n+1}}{2 \left(n+1\right)},
\quad
-\frac{-2 n^{2}+\sqrt{13 n^{2}-22 n+1}+5 n+1}{2 \left(n+1\right)}.
$$
It can be verified that 
$$
(n-1)-\frac 12>
\frac{2 n^{2}-5 n-1+\sqrt{13 n^{2}-22 n+1}}{2 \left(n+1\right)}
>
-\frac{-2 n^{2}+\sqrt{13 n^{2}-22 n+1}+5 n+1}{2 \left(n+1\right)}.
$$
Hence, by \eqref{se4-e2}, 
for $n\ge 3$ and $x\in (n-1-1/2,n-1)$, 
$$
P_{G_{n-3,n}}(x)^2 
<
P_{G_{n-3,n}}(x+1)
P_{G_{n-3,n}}(x-1).
$$
}
\fi

\section{Further problems}
\label{further}

Let $G$ be a graph of order $n\ge2$, and let $\omega(G)$ and $\chi(G)$
denote its clique number and chromatic number, respectively. By Theorem~\ref{main},
\begin{equation}
	\label{e:known}
	P_G(x)^2\ge P_G(x-1)P_G(x+1)
	\qquad\text{for every real }x\ge n-1,
\end{equation}
where the bound $n-1$ is best possible for $n\ge 2$ (see \eqref{complete}).

It is nevertheless natural to seek smaller bounds for suitable classes
of graphs. 
A replacement of $n-1$ by $\omega(G)$ is not valid in general.
Seymour~\cite{Seymour} constructed graphs $H$ such that
$$
	P_H(6)^2<P_H(5)P_H(7).
$$
These examples have $\chi(H)=3$; see also
\cite[Section~1]{Fadnavis}. Hence $6\ge\omega(H)$, and so the
inequality in \eqref{e:known} need not hold for all real
$x\ge\omega(G)$. Theorem~\ref{main} also implies that these counterexamples have
order at least $8$. Consequently, restricting the integer values to
$k\in[\omega(G),n]$ does not avoid these counterexamples.

For chordal graphs, however, the proposed clique-number bound does hold.
To see this, take a perfect elimination ordering of any chordal graph $G$ and let $d_i$ be the
number of neighbors of its $i$-th vertex occurring later in that ordering.
The standard factorization of $P_G(x)$ is
$$
	P_G(x)=\prod_{i=1}^n(x-d_i),
	\qquad \text{where}~0\le d_i\le\omega(G)-1.
$$
%see \cite[Section~5, Proposition~12]{OhPostnikovYoo2008}.
For $x\ge\omega(G)$, all $x-d_i\ge1$, and therefore
$$
	\begin{aligned}
		P_G(x-1)P_G(x+1)
		&=\prod_{i=1}^n\bigl((x-d_i)^2-1\bigr)
		<\prod_{i=1}^n(x-d_i)^2
		=P_G(x)^2.
	\end{aligned}
$$
This suggests the following problem.

\begin{problem}\label{prob:real-clique-bound}
	Identify classes of graphs, beyond chordal graphs, for which
$$
		P_G(x)^2\ge P_G(x-1)P_G(x+1)
		\qquad\text{for every real }x\ge\omega(G).
$$
\end{problem}

For the corresponding discrete problem, Theorem~\ref{main} covers all integers
$k\ge n-1$. Moreover, if $1\le k\le\chi(G)$, then $P_G(k-1)=0$, so the
inequality holds automatically. Thus only the integers satisfying
$\chi(G)+1\le k\le n-2$ require further consideration. Since Seymour's
examples show that the inequality can fail in this range, the remaining
question is to characterize the graphs for which it holds.

\begin{problem}\label{prob:integer-bound}
	Characterize the graphs $G$ of order $n\ge2$ for which
$$
		P_G(k)^2\ge P_G(k-1)P_G(k+1)
$$
	holds for every integer $k$ satisfying
	$\chi(G)+1\le k\le n-2$.
\end{problem}

\medskip
\noindent\textbf{Declaration on the use of AI:}
During the preparation of this work, the authors used AI systems to assist in exploring potential approaches to the proof of Theorem~\ref{main}. All AI-generated suggestions were
verified and refined by the authors, who take full responsibility for the correctness and originality of the paper.

\noindent {\bf Data availability statement}:  Not applicable.

\section*{Acknowledgment}
This work is supported by the Fundamental Research Funds for the Central Universities (No. 20720260153), the National Natural Science Foundation of China (No. 12501498), 
the Natural Science Foundation of Xiamen, China (No. 3502Z202571026), and the Foundation for Cultivated Young Talents of Fujian Province, China (No. 2025350064).

\end{document}